\documentclass{IEEEtran}
\usepackage{cite}
\usepackage{amsmath,amssymb,amsfonts}
\usepackage{graphicx}
\usepackage{textcomp}
\usepackage{bm}
\usepackage{algorithm,algpseudocode}
\usepackage{algorithmicx}
\usepackage{subcaption}

\algrenewcommand\alglinenumber[1]{\footnotesize #1:} 

\def\BibTeX{{\rm B\kern-.05em{\sc i\kern-.025em b}\kern-.08em
    T\kern-.1667em\lower.7ex\hbox{E}\kern-.125emX}}
\begin{document}
\title{Unfitted Lattice Green's Function Method for Exterior Scattering in Complex Geometry}
\author{Siyuan Wang and Qing Xia
\thanks{This work was funded in part by Natural Science Foundation of China (NSFC Grant No: 12401546) and Wenzhou Kean University (Grant No: ISRG2024003 and KY20250604000452). }
\thanks{Siyuan Wang is with Department of Mathematics, Wenzhou Kean University, Zhejiang, China, 325060. (e-mail: sywang@kean.edu). }
\thanks{Qing Xia is with Department of Mathematics, Wenzhou Kean University, Zhejiang, China, 325060. (e-mail: qxia@kean.edu, corresponding author). }
}

\maketitle

\begin{abstract}
This paper develops a finite-difference analogue of the boundary integral/element method for the numerical solution of two-dimensional exterior scattering from scatterers of arbitrary shapes. The discrete fundamental solution, known as the lattice Green’s function (LGF), for the Helmholtz equation on an infinite lattice is derived and employed to construct boundary algebraic equations through the discrete potentials framework. Unlike the continuous fundamental solution used in boundary integral methods, the LGF introduces no singularity, which simplifies numerical implementation. Boundary conditions are incorporated through local Lagrange interpolation on unfitted cut cells. The resulting method retains key advantages of boundary integral approaches—including dimension reduction and the absence of artificial boundary conditions--while enabling finite differences for complex geometries. Numerical results demonstrate the accuracy and robustness of the method for various scatterers, including circular, triangular, and multiple-body configurations.
\end{abstract}


\section{Introduction}\label{sec:introduction}
\IEEEPARstart{F}{inite} difference methods are widely used in the numerical solution of acoustic and electromagnetic scattering problems. Finite difference methods are simple, effective, and natural for simple geometries such as rectangles and cubes. However, irregular shapes with curved boundaries, corners or wedges pose significant geometric challenges for standard finite-difference schemes. To overcome this geometric difficulty, several strategies have been developed. For instance, overset grids \cite{cp1,cp2} were designed such that curvilinear grids are employed near the curved boundaries and away from the boundary are the standard Cartesian grids. 

Another approach is to modify the finite difference stencil at grid points affected by the curved boundary, which is widely used in the immersed interface method \cite{leveque}. When the stencil is not modified, subgrid techniques such as correction based on local truncation error analysis are used to achieve the desired accuracy, as in the kernel-free boundary integral method \cite{ying}, which employs finite difference solutions in the auxiliary domains to approximate the singular integrals, with no needs to evaluate singular kernels explicitly.

Difference potentials method \cite{rya} is another class of finite difference method that can handle curved boundaries or interfaces. The boundary condition is assumed to be represented by spectral basis functions, which are typically selected as Fourier basis functions globally or Chebyshev basis functions for piecewise smooth geometries \cite{med,tsy,britt}. The grid points are connected with the continuous boundary through Taylor's expansion in the normal direction.

Another challenge of using finite differences lies in infinite domains. The infinite domains have to be truncated where artificial boundary conditions are needed for all finite difference methods mentioned above. In contrast, boundary integral/element methods handle infinite domains intrinsically \cite{colton}. No truncation or special numerical tricks are needed in the far field. The outgoing radiation boundary condition is satisfied automatically.

A natural question to ask is if we can design a finite difference version of boundary integral/element method that can handle unbounded regions naturally. Several approaches have been explored.

There is an extensive literature on the computation of the discrete fundamental solution or the lattice Green's function (LGF) of the difference Laplace equation \cite{watson} and a relatively scarce discussion for the LGF of difference Helmholtz equation on the infinite lattice. Nevertheless, the efforts of applying LGF in the numerical solution of exterior scattering can be found.

The first class of methods is fully lattice-based. The boundary is straight and is assumed to align with the grid lines, i.e., the lattice points. The single and double layer formulations can be constructed similarly to the boundary integral framework. The bulk difference equations are reduced to equivalent boundary algebraic equations defined on the boundary lattice points. Once the density at those boundary lattice points are retrieved via incorporating boundary conditions, exterior or interior approximations are reconstructed with means of convolutions with the corresponding single or double layer kernels. As the discrete fundamental solutions share similar structures with their continuous counterparts, fast summation techniques such as Fast Fourier Transform (FFT) or Fast Multiple Method (FMM) can be employed for accelerations of the convolutions. Examples can be found in \cite{dk,jpp1,jpp2}.

The second approach is built upon the first one, but allows the boundary or interface intersects with grid lines in an arbitrary fashion. The boundary difference method developed in \cite{it} and \cite{oa} is the only work that the authors are aware of that combines the LGF with a curved boundary. The boundary difference method was compared extensively with boundary element formulations and emphasizes the singularity-free nature of the LGF. To handle the curved boundary, the higher order Flexible Local Approximation MEthod (FLAME) is employed, which is built on the explicit knowledge of Trefftz basis functions that satisfies the differential equation analytically or approximately. Consequently, the boundary difference method is only applicable to simple geometries circles, ellipses, or squares, as illustrated in the numerical simulations.

The unfitted method developed in this work originates from the difference potentials method \cite{rya}. The primary difference lies in the choice of discrete Green’s function and the use of an indirect formulation.

The main contribution of this work is a new unfitted lattice Green’s function method for 2D exterior Helmholtz scattering that differs from existing lattice-based BAE and boundary difference methods in several aspects: (i) curved boundaries are handled without modified stencils or Trefftz bases, (ii) no singular integrals appear and all kernels are smooth, (iii) the method incorporates boundary conditions through local Lagrange interpolation on cut cells, and (iv) unbounded domains are treated naturally without artificial boundary conditions.

\section{Lattice Green's Function}

In this section, we give a review on the computation of Helmholtz lattice Green's function.

Consider the Helmholtz equation with wavenumber $k$:
\begin{align}
\Delta u +k^2 u = 0
\end{align}
in the free space. We will discuss how to incorporate boundary conditions in later sections. The finite difference discretization on the infinite lattices with grid spacing $h$ translates to the difference Helmholtz equation on $\mathbb{Z}^2$:
\begin{align}\label{eqn:difference_helmholtz}
[Au](m_1,m_2)=0.
\end{align}
Here, $A$ is the standard 5-point central finite difference stencil:
\begin{align}
[Au](m_1,m_2):=&u_{m_1+1,m_2}+u_{m_1,m_2+1}+u_{m_1-1,m_2}\nonumber\\
&+u_{m_1,m_2-1}+(w^2-4)u_{m_1,m_2},
\end{align}
where $u_{m_1,m_2}\approx u(m_1h,m_2h)$, and $\omega^2=k^2h^2$ is the lattice scaled wavenumber.

The forward Fourier transform of the grid function $u_{i,j}$ on the infinite lattice is:
\begin{align}
    \hat{u}_{\xi_1,\xi_2} = \sum_{m_1=-\infty}^{\infty}\sum_{m_2=-\infty}^{\infty} e^{-i\bm{\xi}\cdot\bm{m}}u_{m_1,m_2},
\end{align}
where $\bm{m}=(m_1,m_2)$ and $\bm{\xi}=(\xi_1,\xi_2)$.
The corresponding inverse Fourier transform is given by
\begin{align}
    u_{m_1,m_2} = \frac{1}{(2\pi)^2}\int_{-\pi}^{\pi}\int_{-\pi}^{\pi}e^{i\bm{\xi}\cdot\bm{m}}\hat{u}(\bm{\xi})\,d\bm{\xi}.
\end{align}

To obtain the fundamental solution $G(m_1,m_2)$, we consider
\begin{align}\label{eqn:delta}
[AG](m_1,m_2) = \delta_0(m_1,m_2).
\end{align}
Without loss of generality, here we let $\delta_0$ be the discrete delta function centered at the origin $(0,0)$. Note that the forward transform gives 
\begin{align}
    \hat{\delta_0}(\xi_1,\xi_2) = \sum_{m_1,m_2=-\infty}^{\infty} e^{-i\bm{\xi}\cdot\bm{m}}\delta_0(m_1,m_2) = 1,
\end{align}
since the discrete delta function is 0 except at the origin.
The inverse transform gives
\begin{align}
\delta(m_1,m_2) = \frac{1}{(2\pi)^2}\int_{-\pi}^{\pi}\int_{-\pi}^{\pi}e^{i\bm{\xi}\cdot\bm{m}}\,d\xi_1d\xi_2.
\end{align} 
In the physical space, \eqref{eqn:delta} becomes
\begin{align}
\frac{1}{(2\pi)^2}\int_{-\pi}^{\pi}\int_{-\pi}^{\pi}\left(2\sum_{i=1}^2\cos(\xi_i)+k^2-4\right)e^{i\bm{\xi}\cdot\bm{m}}\hat{G}(\bm{\xi})\,d\bm{\xi}\nonumber\\
=\frac{1}{(2\pi)^2}\int_{-\pi}^{\pi}\int_{-\pi}^{\pi}e^{i\bm{\xi}\cdot\bm{m}}\,d\bm{\xi}.
\end{align}
In the Fourier space, we have
\begin{align}
\left[2\cos(\xi_1)+2\cos(\xi_2)+k^2-4\right]\hat{G}(\bm{\xi})=1,
\end{align}
hence
\begin{align}
\hat{G}(\bm{\xi})=\frac{1}{2\sum_{i=1}^2\cos(\xi_i)+k^2-4}.
\end{align}

The inverse transform would give the lattice Green's function for the difference Helmholtz equation:
\begin{align}
G(\bm{m})=\frac{1}{(2\pi)^2}\int_{-\pi}^{\pi}\int_{-\pi}^{\pi}\frac{\cos(m_1\xi_1)\cos(m_2\xi_2)}{2\sum_{i=1}^2\cos(\xi_i)+k^2-4}\,d\bm{\xi}.
\end{align}
Other stencils are also studied in \cite{Bam}.

$G(m_1,m_2)$ is an even function in $m_1$ and $m_2$ so that we will only consider $m_1\geq0$ and $m_2\geq0$. Also
\begin{align}
G(m_1,m_2) = G(m_2,m_1),
\end{align}
hence we focus on $m_1\geq m_2\geq0$.

Note that in \cite{martin}, the diagonal elements of the lattice Green's function $G$ can be explicitly evaluted in terms of Legendre functions: for $k^2<4$,
\begin{align}\label{eqn:diag1}
G(n,n) = \frac{(-1)^n}{2\pi i}Q_{n-1/2}(z)-\frac{\pi i}{2}P_{n-1/2}(z),
\end{align}
and for $k^2>4$,
\begin{align}\label{eqn:diag2}
G(n,n) = \frac{(-1)^n}{2\pi i}Q_{n-1/2}(z)+\frac{\pi i}{2}P_{n-1/2}(z),
\end{align}
where $z=1-(4-k^2)^2/8$. Note that $Q_{n-1/2}(z)$ blows up at $z=1$, so the above expression is invalid for $k^2=4$, but it can be checked, also noted in \cite{martin} and \cite{bhat}, that a particular solution with alternating rings of $\pm1/4$ exists for $k^2=4$. Nevertheless, $k^2=4$ will not be problematic as we can always nudge away from 4 by choosing a different $h$. The far field asymptotics were studied in \cite{martin}.

It is worthy to note that in \cite{morita}, recursion formulas are given based on the difference equation~\eqref{eqn:difference_helmholtz} and the symmetry property of $G(i,j)$, assuming $G(n,n)$ known.
\begin{align}
G(1,0) = \frac{1}{4}-\frac{1}{4}(k^2-4)G(0,0).
\end{align}
For $j=0$,
\begin{align}
G(m,0) = &(4-k^2)G(m-1,0)\nonumber\\
&-G(m-2,0)-2G(m-1,1).
\end{align}
For $0<n<m-1$,
\begin{align}
G(m,n) = &(4-k^2)G(m-1,n)-G(m-2,n)\nonumber\\
&-G(m-1,n+1)-G(m-1,n-1).
\end{align}
For $n=m-1$,
\begin{align}
G(m,n) &= \frac{4-k^2}{2}G(m-1,n)-G(m-1,n-1).
\end{align}
Unfortunately the recursion formulas are numerically unstable unless high precision floating-point number is used.

To speed up the calculation of the Helmholtz lattice Green's functions, we follow \cite{it} and \cite{oa} and use the continuous Helmholtz Green's function as the boundary condition. The exact procedure is to discretize a box $[-Nh,Nh]\times[-Nh,Nh]$ with same lattices spacing $h$ and for sufficiently large $N$. We approximate the Helmholtz LGF by solving 
\begin{subequations}\label{eqn:LGF}
\begin{align}
[Au](m_1,m_2) &= \delta_{m_1,m_2}, \\
u_{m_1,\pm(N+1)} &= -\frac{i}{4}H^{(1)}_0(\omega h\sqrt{(N+1)^2+m_1^2}),\\
u_{\pm(N+1),m_2} &= -\frac{i}{4}H^{(1)}_0(\omega h\sqrt{(N+1)^2+m_2^2}),
\end{align}
\end{subequations}
where $-N\leq m_1,m_2\leq N$, and $\delta_{m_1,m_2}$ is 1 at the origin and 0 elsewhere. Here the fundamental solution is taken to be the out-going fundamental solution, in contrast to the in-coming one in \cite{it}. The discretizations \eqref{eqn:LGF} can be solved using fast sine transforms. The accuracy testing against \eqref{eqn:diag1} or \eqref{eqn:diag2} when $N>200$ is about $10^{-6}$, which is sufficient given phase error might be dominating for large wave number. Increasing $N$ does not further improve the accuracy though.

\section{Numerical Scheme}

\subsection{Point sets}
The numerical method developed in this work is essentially a finite difference method. To formulate the boundary algebraic equations, we need to first introduce a few point sets. Given a bounded domain $\Omega\subset\mathbb{R}^2$, we embed it into an infinite lattice $(\mathbb{Z}h)^2$ where $h$ is the uniform grid spacing. At each point $\bm{x}_{mn}=(x_m,y_n)$ on the lattice, we can define a 5-point stencil
\begin{align}
N_{mn} = \{(x_m,y_n), (x_{m\pm1},y_n), (x_m,y_{n\pm1})\}.
\end{align}
We will focus on the second order scheme based on the 5-point stencil. High order schemes can be constructed similarly once the corresponding lattice Green's functions are obtained.

\begin{figure}[htbp]
\centering
\includegraphics[width=0.3\textwidth]{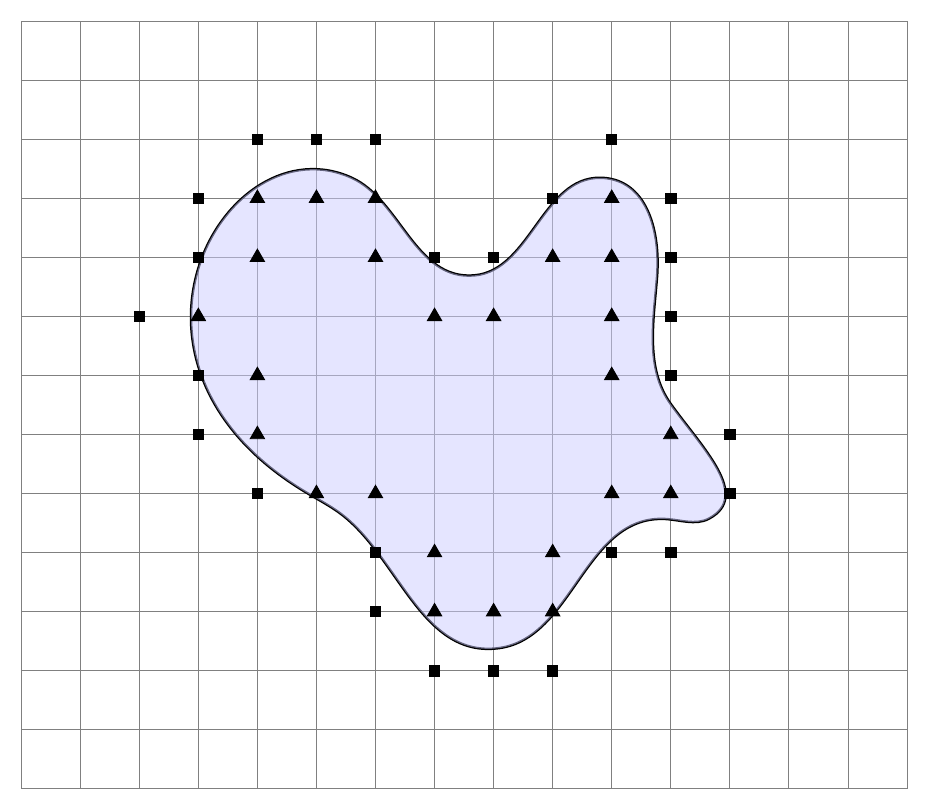}
\caption{Example of set $\gamma$ ($\gamma_{-}$: triangles; $\gamma_{+}$: squares)}\label{fig:gamma}
\end{figure}

The following sets can be defined
\begin{subequations}
\begin{align}
M^+ &=\{(x_m,y_n)\,\vert\,(x_m,y_n)\in\Omega^c\},\\
M^- &=\{(x_m,y_n)\,\vert\,(x_m,y_n)\in\Omega\},\\
N^\pm &=\{N_{mn}\,\vert\,(x_m,y_n)\in M^\pm\},\\
\gamma &= N^+\cap N^-.
\end{align}
\end{subequations}
Here $\Omega^c$ denotes the complement of $\Omega$. In particular, we differentiate between the interior ($-$) and exterior ($+$) parts of the discrete grid boundary $\gamma_{\pm} = \gamma\cap M^\pm$. Fig.~\ref{fig:gamma} shows an example of such points.

Fig.~\ref{fig:gamma} shows that the point set $\gamma$ straddles the continuous boundary. In addition, an extra interior layer can be defined: 
\begin{align}
E = \{N_{mn}\,|\,(x_m,y_n)\in \gamma_-\}\,\backslash\, \gamma,
\end{align}
which includes additional interior points that are connected to points in $\gamma_-$ from the 5-point stencil.

\subsection{Indirect formulations} In this part, we mimic the boundary integral method and formulate the boundary algebraic version. Similar formulations can be found in \cite{tsyn2,lgf,zheng}.
\paragraph{Single-layer formulation} The single layer kernel is defined as
\begin{align}
S(s,t) = G(s-t),
\end{align}
for $s,t\in\mathbb{Z}^2$. For interior points in $N^+$, we have the following single layer representation
\begin{align}\label{eqn:single}
u(s) = \sum_{\bm{x}_t\in \gamma_-}S(s,t)q_t,\quad \bm{x}_s\in N^+,
\end{align}
where $q_t$ is the unknown single-layer densities defined at $\gamma_-$. It can be shown that the single layer formulation  satisfies the homogeneous difference Helmholtz equation.

\paragraph{Double-layer formulation} The double layer kernel is defined as
\begin{align}
D(s,t) = \sum_{k\in E_t}G(s-t)-G(s-r),
\end{align}
for $s,t,r\in\mathbb{Z}^2$. $E_t$ denotes the points in $E$ and connected to the point $t$. For interior points in $N^+$, we have the following single layer representation
\begin{align}\label{eqn:double}
u(s) = \sum_{\bm{x}_t\in \gamma_-}D(s,t)q_t,\quad \bm{x}_s\in N^+,
\end{align}
where $q_t$ denotes the unknown double-layer densities defined at $\gamma_-$. It also can be shown that the double layer formulation  satisfies the homogeneous difference Helmholtz equation.

As studied in \cite{jpp1}, a combined kernel in the form of
\begin{align}
C(s,t) = D(s,t)-i\eta S(s,t),
\end{align}
which is a discrete analogue of Burton-Miller combined integral,
can be constructed to suppress spurious oscillations. The parameter $\eta$ will play to minimize the condition number of the coefficient matrix. We will not explore the combined form in this work.

Given a kernel $K(s,t)$, which can be $S(s,t),D(s,t)$ or $C(s,t)$, we can relate $\gamma_\pm$ via the unknown density defined at $\gamma_-$:
\begin{subequations}
\begin{align}
u_{\gamma_-} = K_-q,\\
u_{\gamma_+} = K_+q,
\end{align}
\end{subequations}
where $K_\pm$ are the coefficient matrix for the convolutions. $K_-$ is a square, dense, and invertible matrix, which leads to
\begin{align}\label{eqn:relation}
u_{\gamma_+} = K_+K_-^{-1}u_{\gamma_-}.
\end{align}
The linear relation can be used to design artificial boundary conditions \cite{yin}. We note that $K_-$ is dense and inverting dense matrix is unfeasible for large degree of freedoms.

\subsection{Boundary closure}

The relation in \eqref{eqn:relation} only exists on the discrete grid boundary $\gamma$ and is equivalent to the Helmholtz difference equation. In this part, we will discuss how to close the system with the boundary conditions.

\begin{figure}[htbp]
\centering
\includegraphics[width=0.2\textwidth]{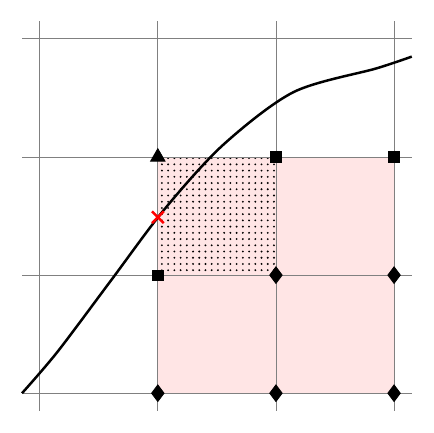}
\caption{An example of a cut cell ($\gamma_{-}$: triangles; $\gamma_{+}$: squares; $\zeta$: diamonds)}\label{fig:cut}
\end{figure}

To this end, we introduce basis functions similarly to the finite element method. For each point (the triangle point in Fig.~\ref{fig:cut}) in $\gamma_-$, we identify an associated intersection point (the cross point in Fig.~\ref{fig:cut}) between the grid line and the boundary $\Gamma$. 

Then, based on the dotted cut cell, we locate a local $3\times3$ mesh (see Fig.~\ref{fig:cut} for an example). We use as many points in the computational domain as possible when selecting the extra grid points $\zeta$ (diamond points). The point set $\gamma_+$ is then augmented with $\zeta$ and we define $\tilde{\gamma}_+:=\gamma_+\cup\zeta$ and $\tilde{\gamma}:=\tilde{\gamma}_+\cup{\gamma}_-$.

On the local mesh which consists of points in $\gamma_-$ and $\tilde{\gamma}_+$, we define local Lagrange basis functions
\begin{align}
\phi_{mn}(x,y) = \phi_m(x)\phi_n(y),\quad m,n=1,2,3,
\end{align}
where the barycentric form of the 1D Lagrange basis function is
\begin{align}
\phi_m(x) = \ell(x)\frac{\omega_m}{x-x_m},\quad m=1,2,3,
\end{align}
where
\begin{align}
\ell(x) = \prod_{m=1}^3(x-x_m),\quad \omega_m= \prod_{n\neq m} \frac{1}{x_n-x_m}=\frac{1}{h^2},
\end{align}
where $h$ is the uniform grid size. In $y$ direction, $\phi_n(y)$ is constructed similarly.

We assume the boundary conditions are in the form of $\alpha u_{\bm{n}}+\beta u = g$, and  $\bm{n}$ is the outward unit normal vector. The discretization would lead to
\begin{align}
\sum_{x_{mn}\in \tilde{\gamma}} u_{mn}\left[\alpha \nabla \phi_{mn}(x) \cdot \bm{n}(x)+\beta \phi_{mn}(x)\right] = g(x),
\end{align}
where $x$ denotes intersection points.
The above form can also be expressed into a matrix-vector form $\Phi u_{\tilde{\gamma}}=g$ where $\Phi$ is a sparse matrix. It can also be decomposed into the following form
\begin{align}
\Phi_+u_{\tilde{\gamma}_+}+ \Phi_{-}u_{{\gamma}_-} = g(x).
\end{align}
Since we add extra points in $\zeta$, we modify \eqref{eqn:relation}, combine the boundary condition, and we reach
\begin{subequations}\label{eqn:boundary_system}
\begin{align}
u_{\tilde{\gamma}_+} = \tilde{K}_+K_-^{-1}u_{\gamma_-},\\
\Phi_+u_{\tilde{\gamma}_+}+ \Phi_{-}u_{{\gamma}_-} = g(x).
\end{align}
\end{subequations}
Here, the tilde denotes matrix for the augmented point set.

We take two strategies to solve the boundary algebraic equations \eqref{eqn:boundary_system}. The first is to use Gaussian elimination and eliminate $u_{\gamma_-}$, which is also known as the Schur complement:
\begin{align}
(\Phi_+\tilde{K}_+K_-^{-1}+ \Phi_{-})u_{{\gamma}_-} = g(x).
\end{align}
This approach is more expensive as we need to invert $K_-$, but would result in a well-conditioned system similar to the Fredholm integral of the second type.

The second approach is to solve for the unknown density $q:=K^{-1}u_{{\gamma}_-}$ defined on $\gamma_-$ from the system
\begin{align}
(\Phi_+\tilde{K}_++ \Phi_{-}K_-)q = g(x).
\end{align}

\subsection{Reconstruction}

In either approach, we aim to solve for the unknown density $q$ defined on $\gamma_-$. Once $q$ is known, the value at any point $s$ can be evaluated using the convolution
\begin{align}
u_s = \sum_{x_t\in\gamma_-}K(s,t)q(t).
\end{align}
For a target region $\omega_h$, it is expensive to evaluate at all points using such convolution. Instead, we can use the convolution to evaluate the values at the boundary
$\partial\omega_h$ only and apply FFT to obtain all interior nodal values in $\omega_h$.

\subsection{Algorithm}

We summarize the flow of the developed method (the second approach) below.
\begin{algorithm}[htbp]
\caption{Unfitted Lattice Green's Function Method}\label{alg:bae}
\allowdisplaybreaks
\begin{algorithmic}[1]

    \State \textbf{Input:} Domain $\Omega$; boundary condition $g$; grid size $h$; wavenumber $k$
    \State \textbf{Output:} Approximate solution $u_s$ in the target region $\omega_h$

    \Statex
    \State $\text{LGF} \gets \text{ComputeLGF}(k,h)$ 
    \Comment{Compute or retrieve the lattice Green's function}

    \State $M^\pm, N^\pm, \gamma^\pm, \tilde{\gamma}_+, x 
        \gets \text{ClassifyPoints}(\Omega,h)$
    \Comment{Identify interior/exterior nodes, cut cells, and intersection points}

    \State $\Phi_\pm \gets \text{ConstructBasisMatrices}(\tilde{\gamma}_+, \gamma_-, x)$
    \Comment{Assemble local basis operators on the augmented stencil}

    \State $\tilde{K}_+, K_- \gets \text{ConstructKernels}(\text{LGF}, \tilde{\gamma}_+, \gamma_-)$
    \Comment{Assemble discrete single-, double-, or combined-layer kernels}

    \State $A \gets \Phi_+ \tilde{K}_+ + \Phi_- K_-$
    \Comment{Form the boundary algebraic system}

    \State $q \gets \text{Solve}(A, g)$
    \Comment{Solve for layer density using a Krylov method (matrix-free)}

    \State $u_s \gets \text{Convolve}(K_{\omega_h}, q)$
    \Comment{Reconstruct the field in $\omega_h$ using FFT-based convolution}

\end{algorithmic}
\end{algorithm}

\subsection{Computational Cost}

The proposed algorithm combines elements of finite element discretization, finite difference schemes, and boundary integral formulations. Step~7 may be interpreted as evaluating a discrete analogue of the boundary integral operator: the action of $K_{-}$ and $\widetilde{K}_{+}$ on the density $q$ provides a second–order approximation of the continuous potentials.

In Step~8, we can solve the linear system $Aq = g$ using a Krylov subspace method. In all numerical experiments, GMRES proved effective and no preconditioner was required. Matrix-free implementation can also be employed in the iterative solver.

The entries of $K_{\pm}$ are obtained directly from the lattice Green's function. The LGF can either be precomputed and stored, when memory permits, or evaluated on the fly during each matrix–vector multiplication. In the latter case, the cost of applying $K_{\pm}$ is $\mathcal{O}(N^{2})$, where $N$ denotes the number of lattice points along one dimension of the computational box. The reconstruction of the field in the bulk region employs FFT-based convolutions and therefore requires $\mathcal{O}(N^{2}\log N)$ operations.

\subsection{Conditioning and Stability}

The stability of the unfitted boundary algebraic formulation depends primarily on the conditioning of the operator $K_{-}$. Numerical experiments indicate that the system matrix $A = \Phi_{+}\tilde{K}_{+} + \Phi_{-}K_{-}$ is moderately well conditioned, and GMRES converges rapidly without preconditioning. This behavior is consistent with classical boundary integral formulations.

As with all finite-difference discretizations of the Helmholtz equation, the choice of grid spacing $h$ must properly resolve the wavelength $\lambda = 2\pi/k$ to avoid numerical dispersion and pollution errors. A commonly accepted requirement is that the grid should provide at least $8$–$12$ points per wavelength. The unfitted LGF-based formulation does not eliminate this resolution constraint, because the boundary algebraic system relies on the standard five-point interior stencil.

\subsection{Memory Usage}

The memory footprint of the method is modest. When precomputation is employed, the storage of the LGF requires $\mathcal{O}(N^{2})$ memory, where $N$ is the number of lattice nodes along one dimension of the truncated computational box. In the matrix‐free implementation, neither $K_{-}$ nor $\tilde{K}_{+}$ is stored explicitly; only point coordinates are retained and all kernel evaluations are performed on demand. Consequently, the memory usage is dominated by the temporary arrays required for FFT-based reconstruction, which scale as $\mathcal{O}(N^{2})$, making the proposed approach suitable for large‐scale simulations.

\section{Numerical Simulations}
In this section, we test the algorithms with the same wave number $k=10$ and planewave $u^{\rm inc}=e^{ik(d_1 x+d_2 y)}$ in the direction of $d=(1/2,\sqrt{3}/2)$. All field plots are results from the single layer formulations as both formulation give identical results. The only difference among these numerical tests is the scatterer geometry.
\subsection{Planewave scattering outside a circle}
First of all, we test with a circular case where analytic solutions can be obtained.
Let the incident direction be the unit vector $d=(\cos\theta_{\rm inc},\sin\theta_{\rm inc})$. In polar coordinates $(r,\theta)$ the plane incidence wave is
\begin{align}
u^{\rm inc}(x,y)=e^{ik r\cos(\theta-\theta_{\rm inc})},
\end{align}
which has the expansion
\begin{align}
u^{\rm inc}(r,\theta)=\sum_{m=-\infty}^{\infty} i^{m} J_{m}(kr)e^{-i m\theta_{\rm inc}}e^{i m\theta},
\end{align}
and $J_m$ is the Bessel function of the first kind.

We seek the scattered field in the form
\begin{align}
u^{\rm sc}(r,\theta)=\sum_{m=-\infty}^{\infty} a_m H^{(1)}_{m}(kr)e^{i m\theta},
\end{align}
where $H^{(1)}_m$ is the outgoing Hankel function. The total field for $r\ge R$ is $u^{\rm tot}=u^{\rm inc}+u^{\rm sc}$.

We will study two variations of Perfect Electric Conductor (PEC) boundary condition for scalar Helmholtz equations: Transverse Magnetic (TM) mode with Dirichlet BC and Transverse Electric (TE) mode with Neumann BC.

\begin{figure}[htbp]
    \centering
    \begin{subfigure}{0.23\textwidth}
        \centering
        \includegraphics[width=\textwidth, trim = 2cm 7cm 2cm 6.5cm]{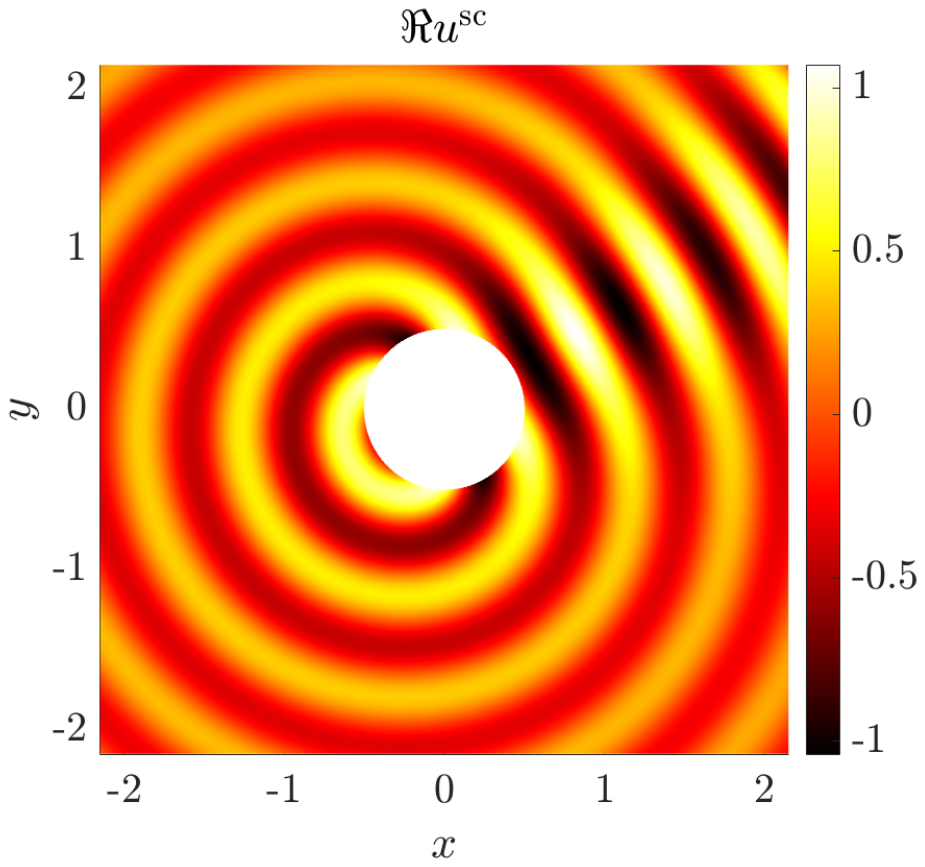}
        \caption{TM: real part of $u^{\rm sc}$}
    \end{subfigure}
    ~
    \begin{subfigure}{0.23\textwidth}
        \centering
        \includegraphics[width=\textwidth, trim = 2cm 7cm 2cm 6.5cm]{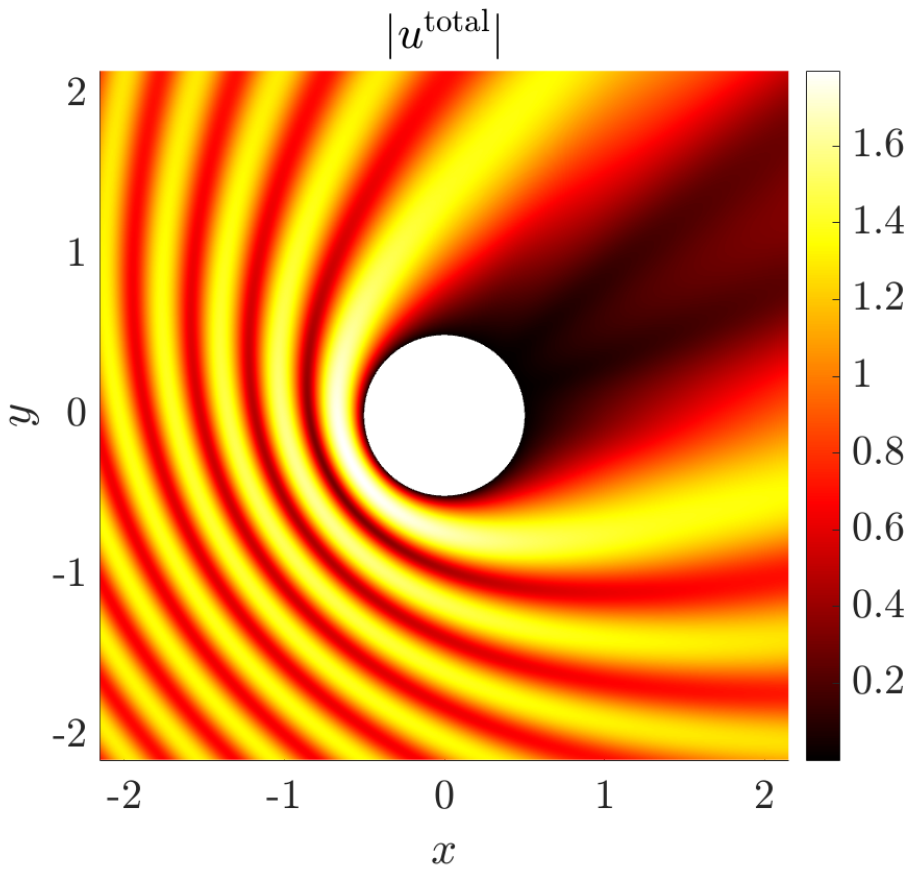}
        \caption{TM: magnitude of $u^{\rm tot}$}
    \end{subfigure}
        \begin{subfigure}{0.23\textwidth}
        \centering
        \includegraphics[width=\textwidth, trim = 2cm 7cm 2cm 6.5cm]{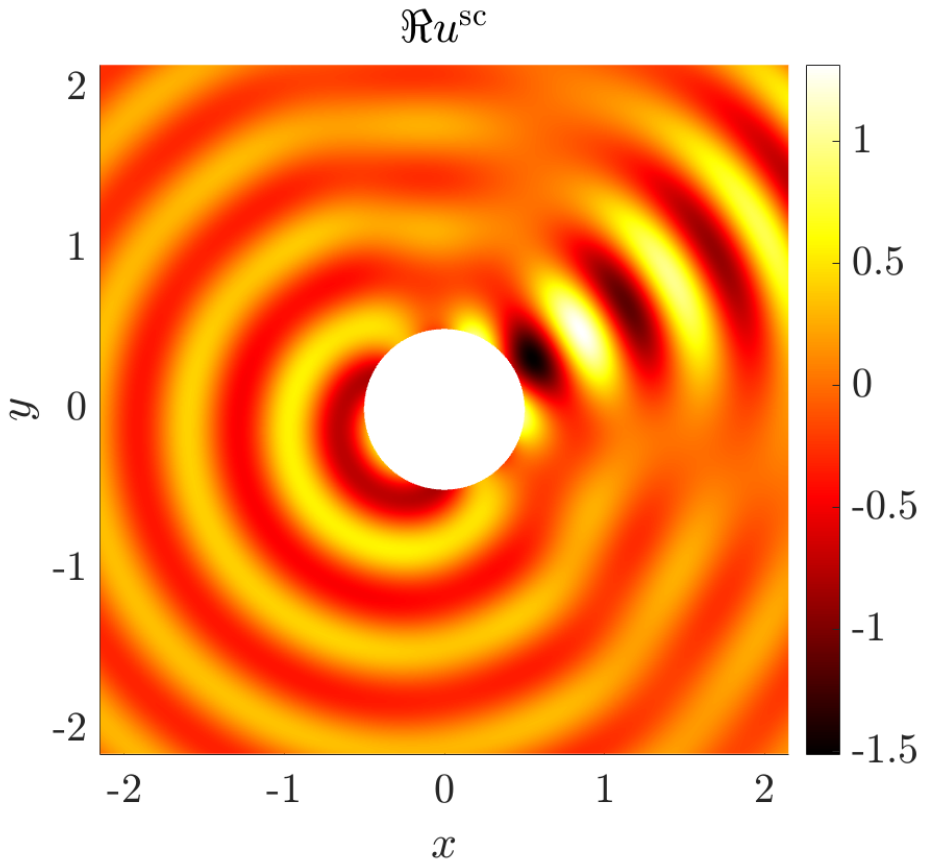}
        \caption{TE: real part of $u^{\rm sc}$}
    \end{subfigure}
    ~
    \begin{subfigure}{0.23\textwidth}
        \centering
        \includegraphics[width=\textwidth, trim = 2cm 7cm 2cm 6.5cm]{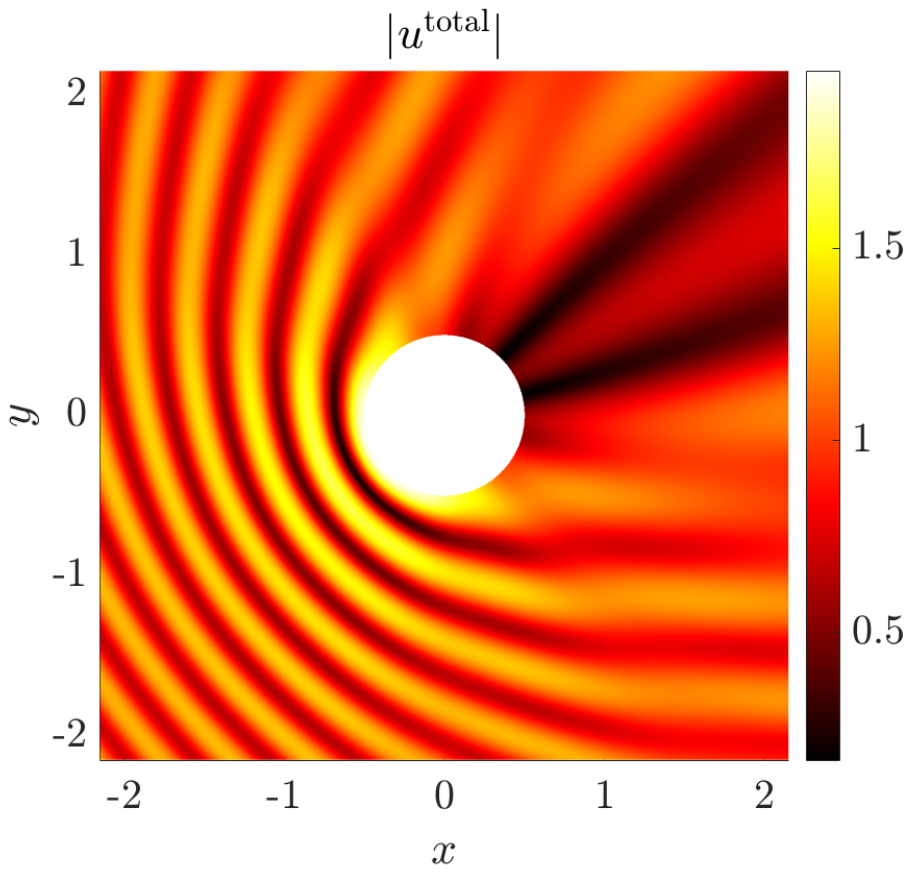}
        \caption{TE: magnitude of $u^{\rm tot}$}
    \end{subfigure}
    \caption{Exact scattering with a circle ($r=0.5,k=10$, single-layer formulation)}\label{fig:planewave}
\end{figure}

\paragraph{TM polarization -- Dirichlet BC}

Boundary condition at $r=R$ gives $u^{\rm tot}(R,\theta)=0$. Enforcing the BC term-by-term yields
\begin{align}
i^{m}J_{m}(kR)e^{-i m\theta_{\rm inc}} + a_m H^{(1)}_{m}(kR)=0,
\end{align}
where we solve for $a_m$ and reach
\begin{align*}
u^{\rm sc}(r,\theta)= -\sum_{m=-\infty}^{\infty} i^{m}\frac{J_{m}(kR)}{H^{(1)}_{m}(kR)}H^{(1)}_{m}(kr)e^{i m(\theta-\theta_{\rm inc})}.
\end{align*}

The series decays rapidly. We truncate it with the largest $|m|=100$ and use the truncation as the ``exact'' solution to compute the errors and convergence.

\paragraph{TE polarization -- Neumann BC}

Boundary condition: normal derivative vanishes on $r=R$, i.e., $\partial_{r}u^{\rm tot}(R,\theta)=0$. Using $\partial_r J_m(kr)=k J'_m(kr)$ and similarly for $H^{(1)}_m$, term-by-term gives
\begin{align}
i^{m}k J'_{m}(kR)e^{-i m\theta_{\rm inc}} + a_mk H^{(1)\prime}_{m}(kR)=0,
\end{align}
where we solve for $a_m$ and reach
\begin{align*}
u^{\rm sc}(r,\theta)= -\sum_{m=-\infty}^{\infty} i^{m}\frac{J'_{m}(kR)}{H^{(1)\prime}_{m}(kR)}H^{(1)}_{m}(kr)e^{i m(\theta-\theta_{\rm inc})},
\end{align*}
which can also be truncated with the largest $|m|=100$.

Fig.~\ref{fig:planewave} shows the field maps of the scattered fields and total fields outside a circle on the finest mesh with $N=2048$. A lens effect can be observed in the scattered field in the TE polarization. The field maps are generated with the single layer formulations, as the other two formulations will give identical results.

\begin{figure}[htbp]
    \centering
    \begin{subfigure}{0.23\textwidth}
        \centering
        \includegraphics[width=\textwidth, trim = 2cm 7cm 2cm 6.5cm]{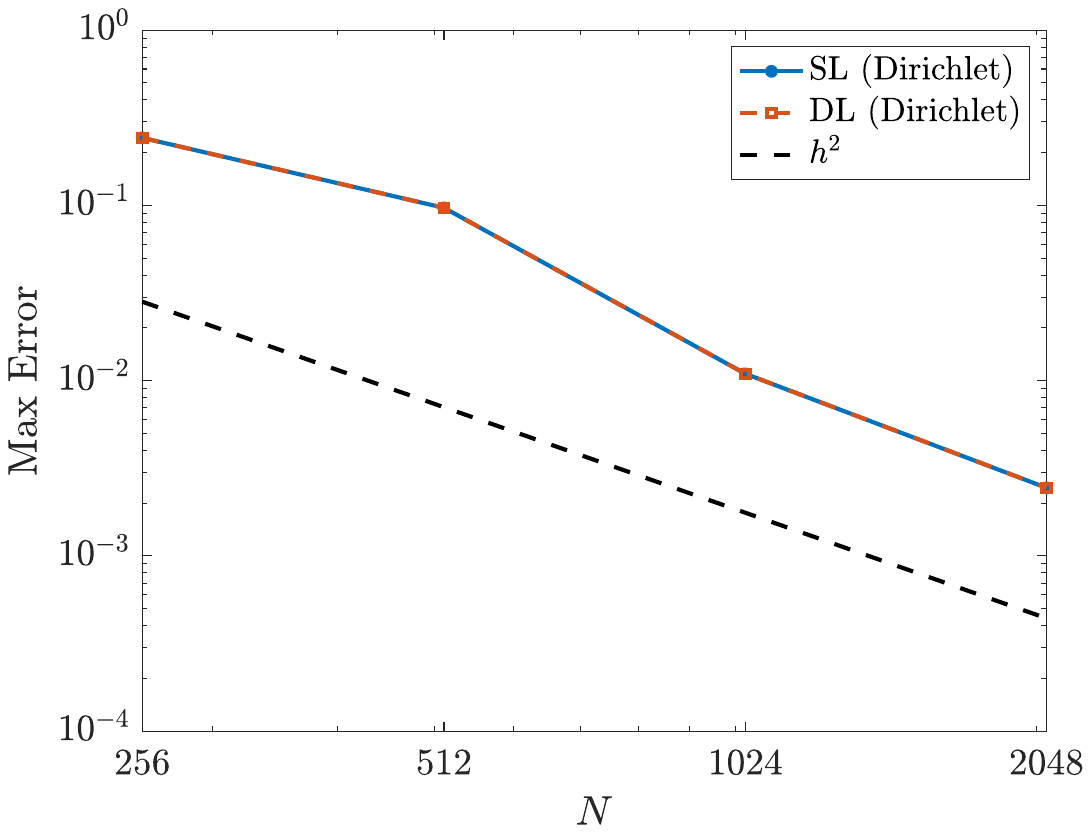}
        \caption{TM}
    \end{subfigure}
    \begin{subfigure}{0.23\textwidth}
        \centering
        \includegraphics[width=\textwidth, trim = 2cm 7cm 2cm 6.5cm]{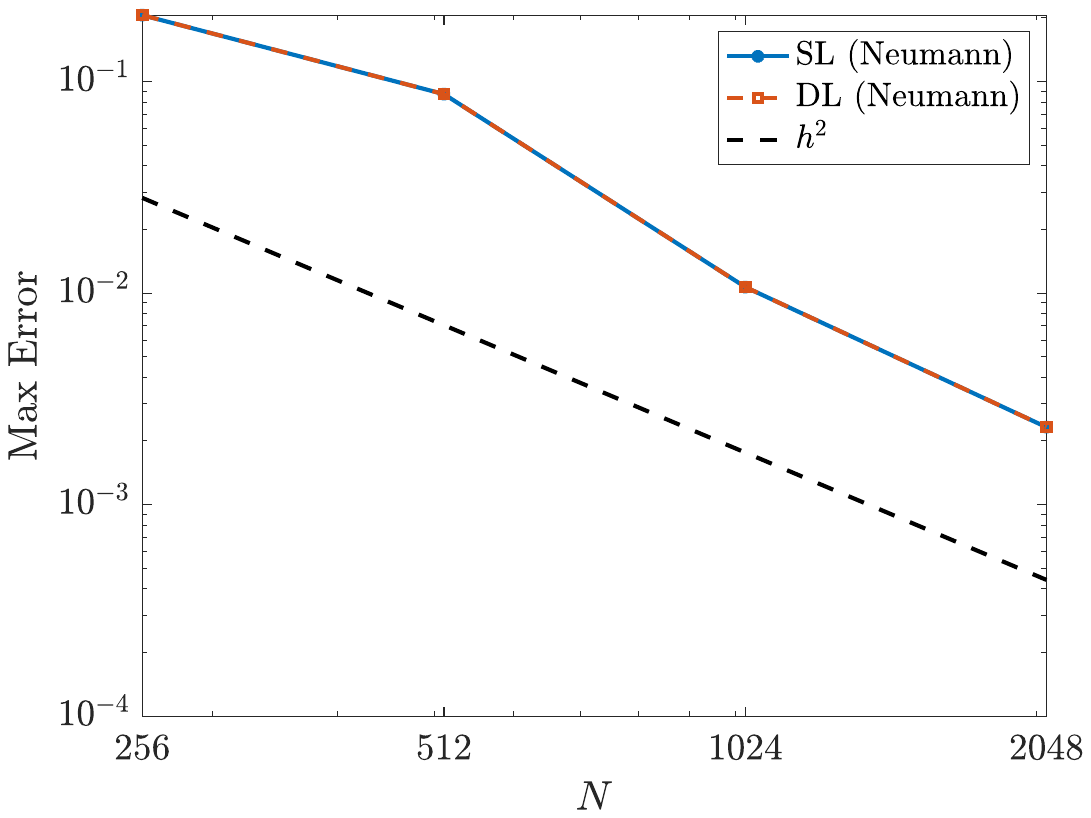}
        \caption{TE}
    \end{subfigure}
    \caption{Convergence ($r=0.5,k=10$); SL: single layer; DL: double layer}\label{fig:conv}
\end{figure}
The convergence in Fig.~\ref{fig:conv} demonstrates second order convergence in max norm in the target domain. The choice of single or double layer formulations does not affect the accuracy. TM and TE polarizations also give similar error magnitude, indicating robustness of the algorithm with respect to different boundary conditions.

\begin{figure}[htbp]
    \centering
    \begin{subfigure}{0.23\textwidth}
        \centering
        \includegraphics[width=\textwidth, trim = 2cm 7cm 2cm 6.5cm]{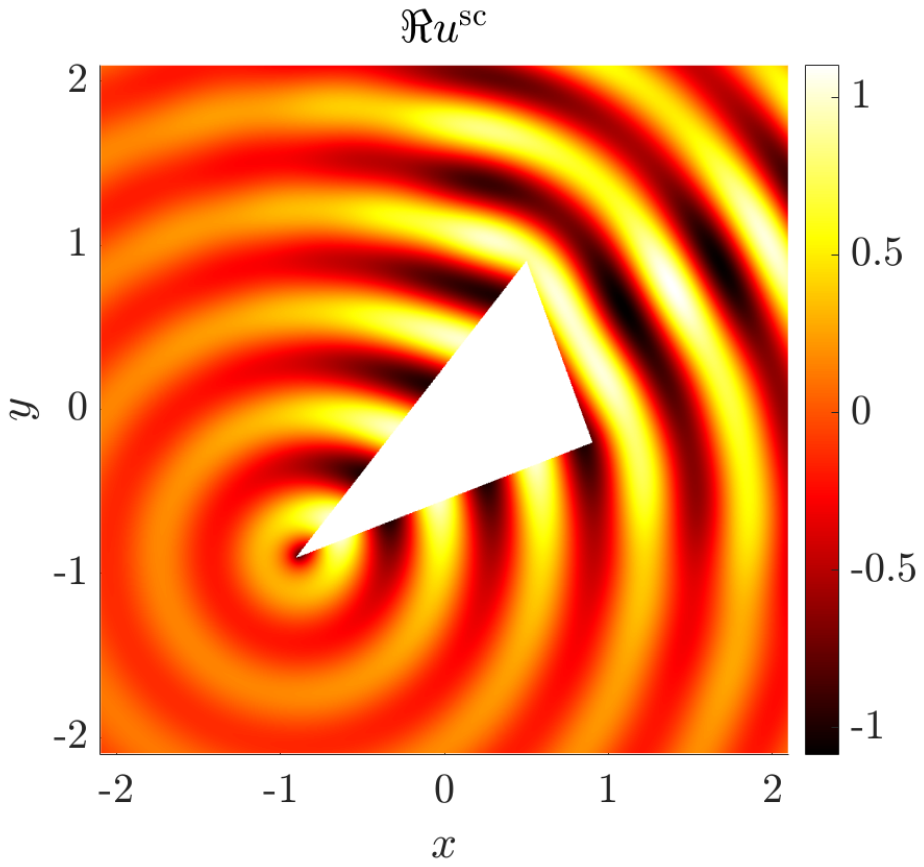}
        \caption{Real part of $u^{\rm sc}$}
    \end{subfigure}
    ~
    \begin{subfigure}{0.23\textwidth}
        \centering
        \includegraphics[width=\textwidth, trim = 2cm 7cm 2cm 6.5cm]{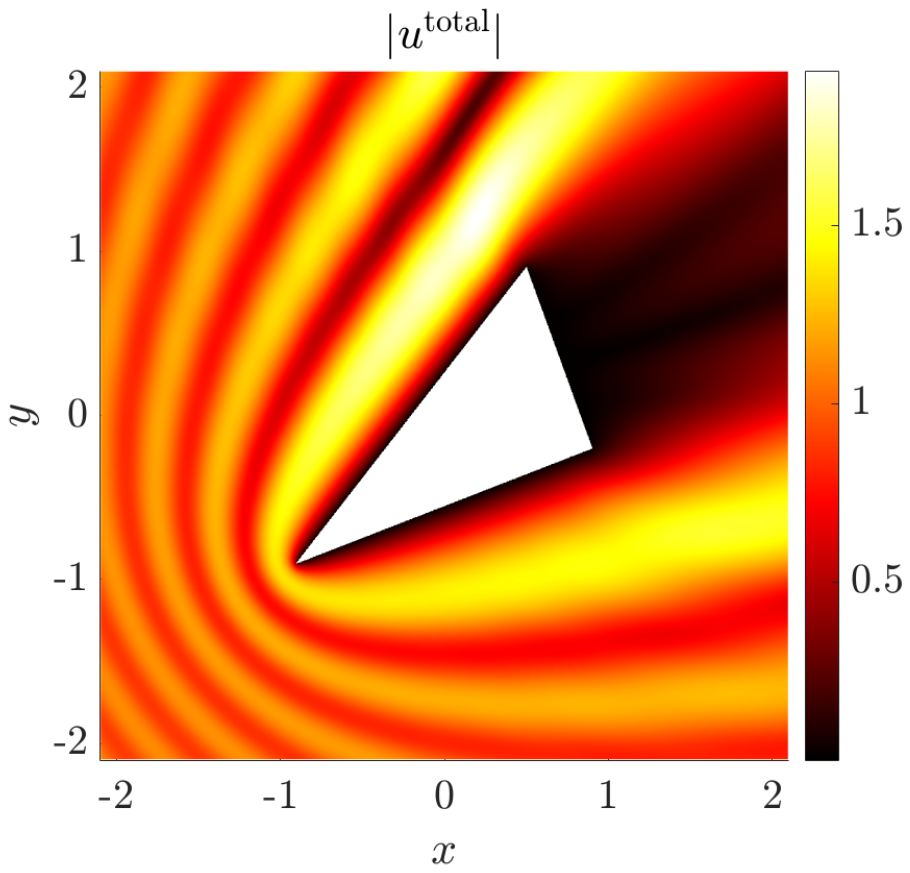}
        \caption{Magnitude of $u^{\rm tot}$}
    \end{subfigure}
    \caption{TM scattering with a triangle}\label{fig:triangle}
\end{figure}

\subsection{Corners}
Next we present a simulation with a triangle, which is challenging for traditional numerical methods. The vertices of the triangle are randomly selected to be $(0.5,0.9), (0.9,-0.2), (-0.9,-0.9)$. No mesh refinement or special numerical technique is needed near the corners. In fact, only point sets $\gamma_\pm$ plus intersection points are needed in the algorithm and the lattice does not ``see'' the corners.

\begin{figure}[htbp]
    \centering
    \includegraphics[width=0.25\textwidth, trim = 2cm 7cm 2cm 6.5cm]{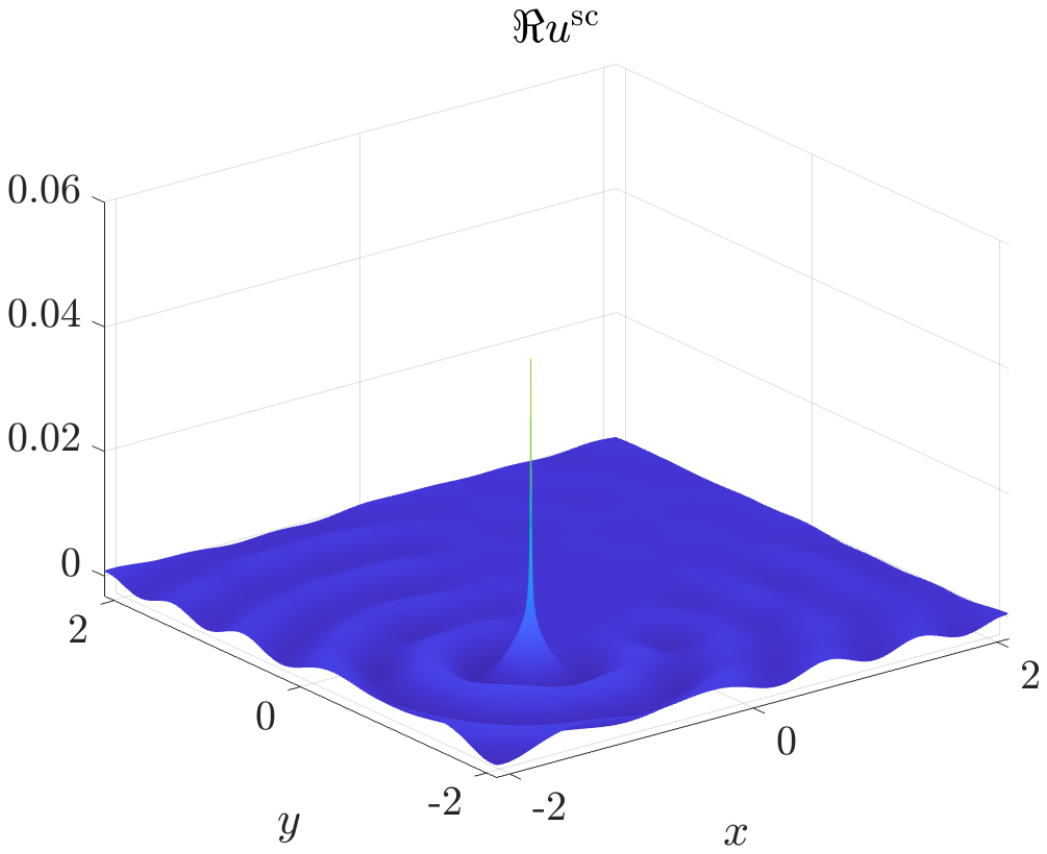}
    \caption{Real part of of the error for $u^{\rm sc}$ for $N=2048$}\label{fig:triangle_error}
\end{figure}

In the scattered field shown in Fig.~\ref{fig:triangle}, a seemingly point source is captured, and the total field is smooth in magnitude. Fig.~\ref{fig:triangle_error} shows the difference of the real part of the scattered field between mesh $N=2048$ and $N=4096$. We observe a sharp single source at a vertex of the triangle, which interrupts the second order convergence.

\subsection{Multiple scattering}

Another merit of the developed algorithm is that it is insensitive to the number of scatterers. The algorithm only needs information of the interior and exterior point sets that are constructed.
The two rods are given by center lines with pair of end points $(-1.5, 0.75)-(-0.25, 0.125)$
and $(0.25, -0.125)-(1.5, -0.75)$, respectively. The radius of the rods is $r = 0.1$.

\begin{figure}[htbp]
    \centering
    \begin{subfigure}{0.23\textwidth}
        \centering
        \includegraphics[width=\textwidth, trim = 2cm 7cm 2cm 6.5cm]{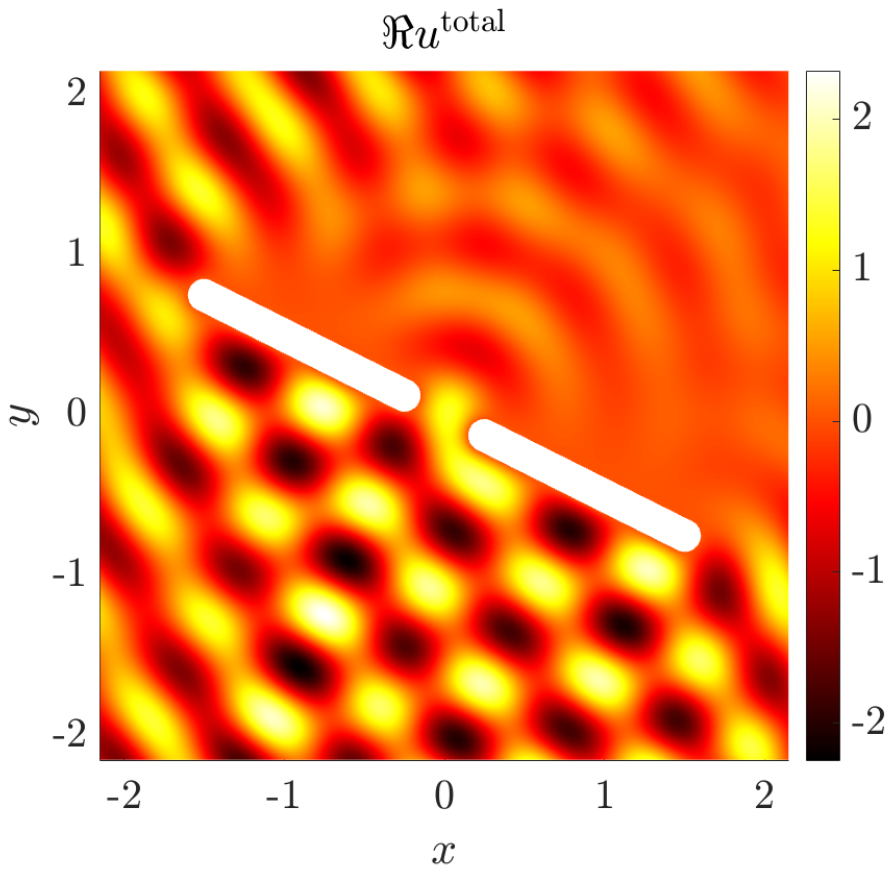}
        \caption{TM: real part of $u^{\rm tot}$}
    \end{subfigure}
    \begin{subfigure}{0.23\textwidth}
        \centering
        \includegraphics[width=\textwidth, trim = 2cm 7cm 2cm 6.5cm]{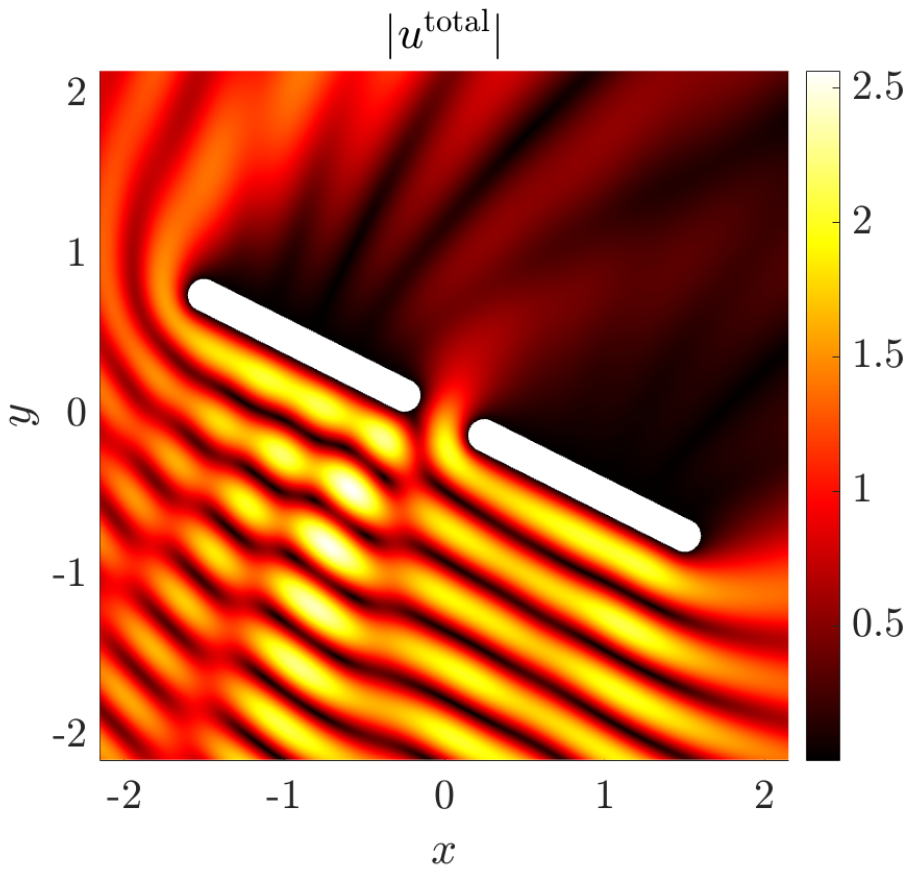}
        \caption{TM: magnitude of $u^{\rm tot}$}
    \end{subfigure}
    \begin{subfigure}{0.23\textwidth}
        \centering
        \includegraphics[width=\textwidth, trim = 2cm 7cm 2cm 6.5cm]{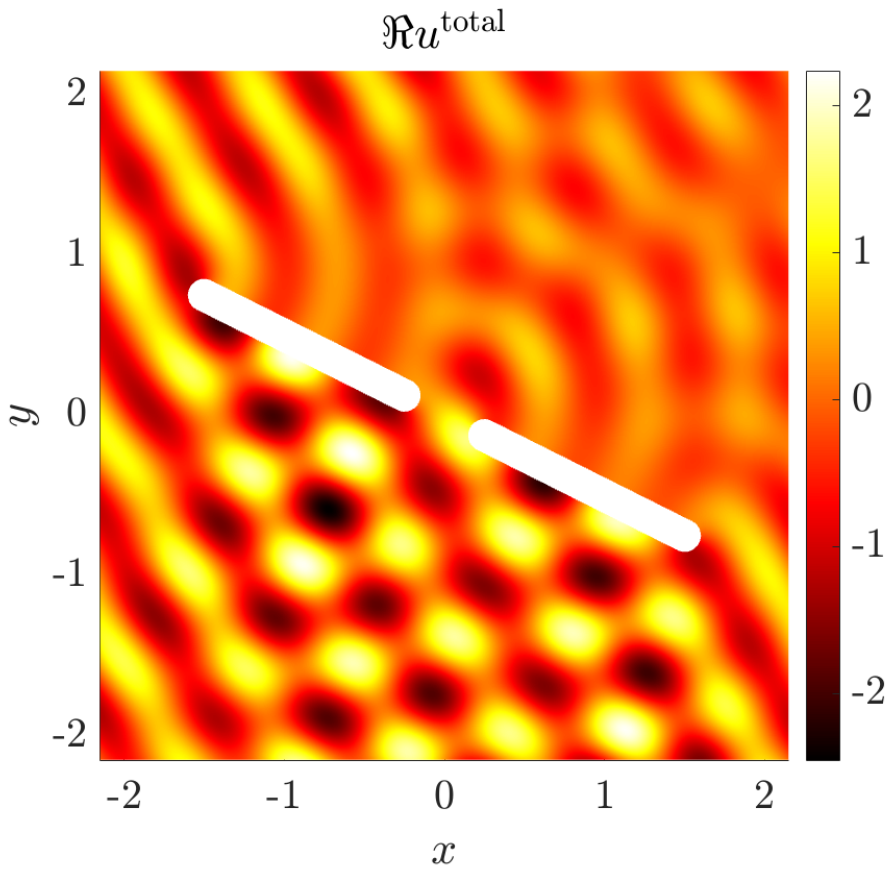}
        \caption{TE: real part of $u^{\rm tot}$}
    \end{subfigure}
    \begin{subfigure}{0.23\textwidth}
        \centering
        \includegraphics[width=\textwidth, trim = 2cm 7cm 2cm 6.5cm]{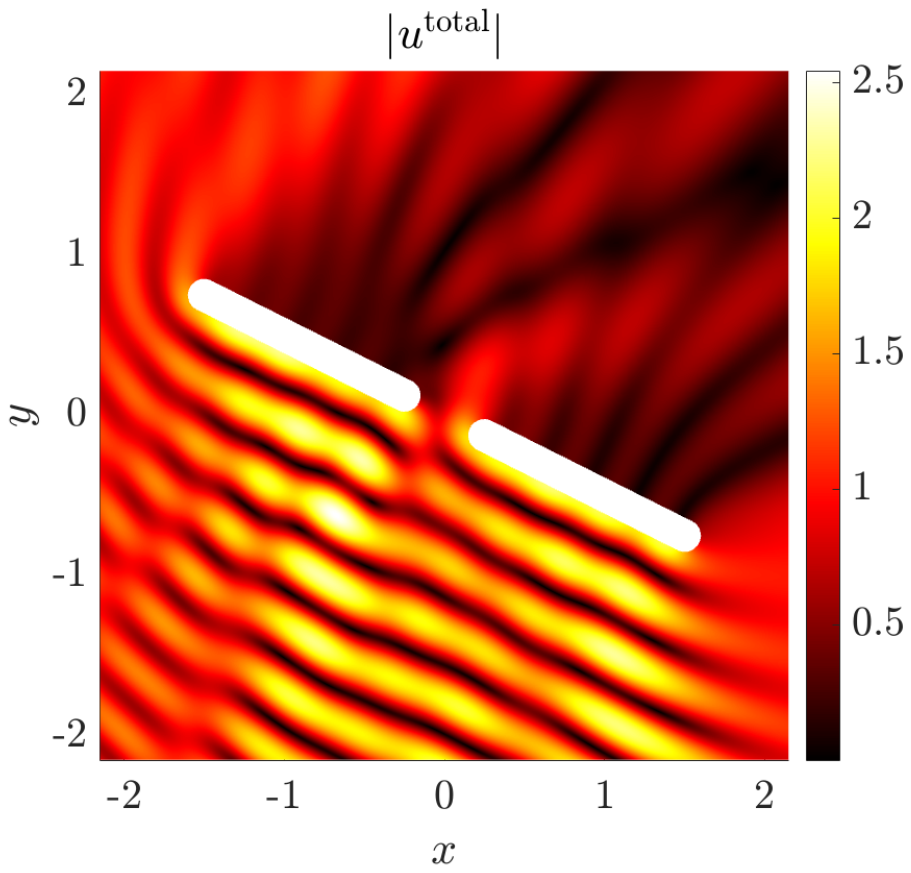}
        \caption{TE: magnitude of $u^{\rm tot}$}
    \end{subfigure}
    \caption{Scattering with two rods}\label{fig:two_rods}
\end{figure}
In both the TM and TE modes, scattering behind the two rods can be observed in the real part of the total field in Fig.~\ref{fig:two_rods}.

\begin{figure}[htbp]
    \centering
    \begin{subfigure}{0.23\textwidth}
        \centering
        \includegraphics[width=\textwidth, trim = 2cm 7cm 2cm 6.5cm]{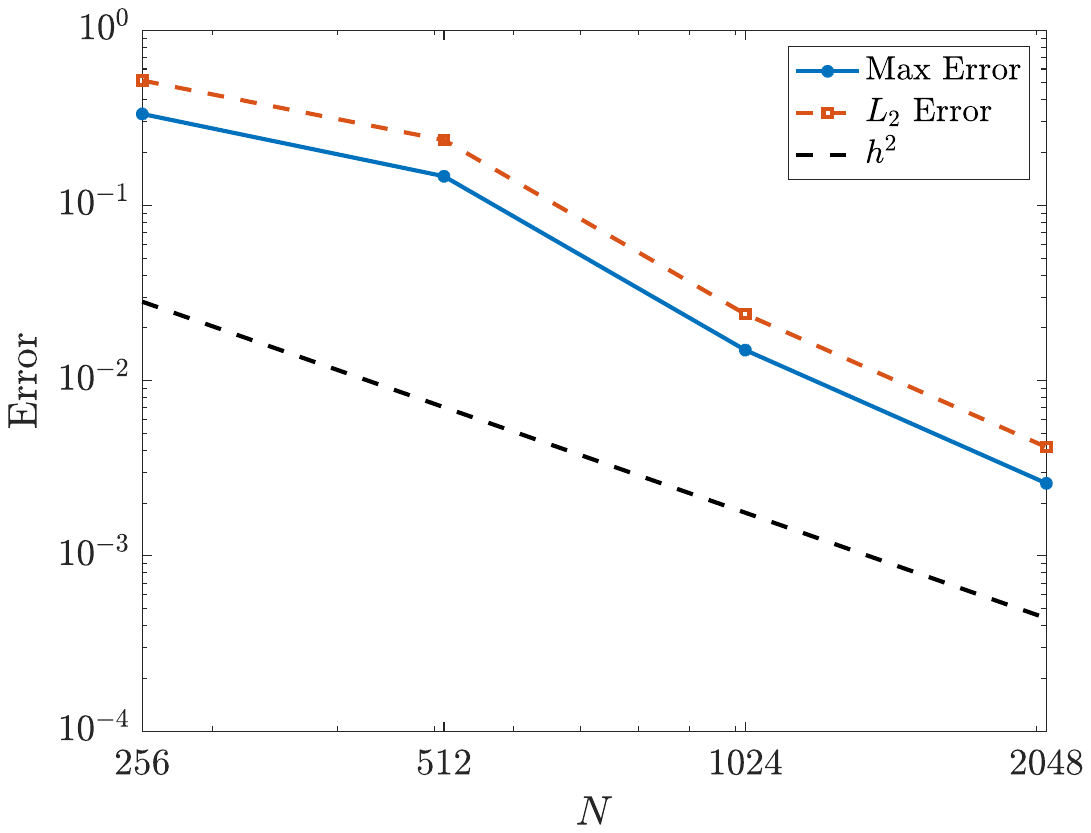}
        \caption{TM}
    \end{subfigure}
    ~
    \begin{subfigure}{0.23\textwidth}
        \centering
        \includegraphics[width=\textwidth, trim = 2cm 7cm 2cm 6.5cm]{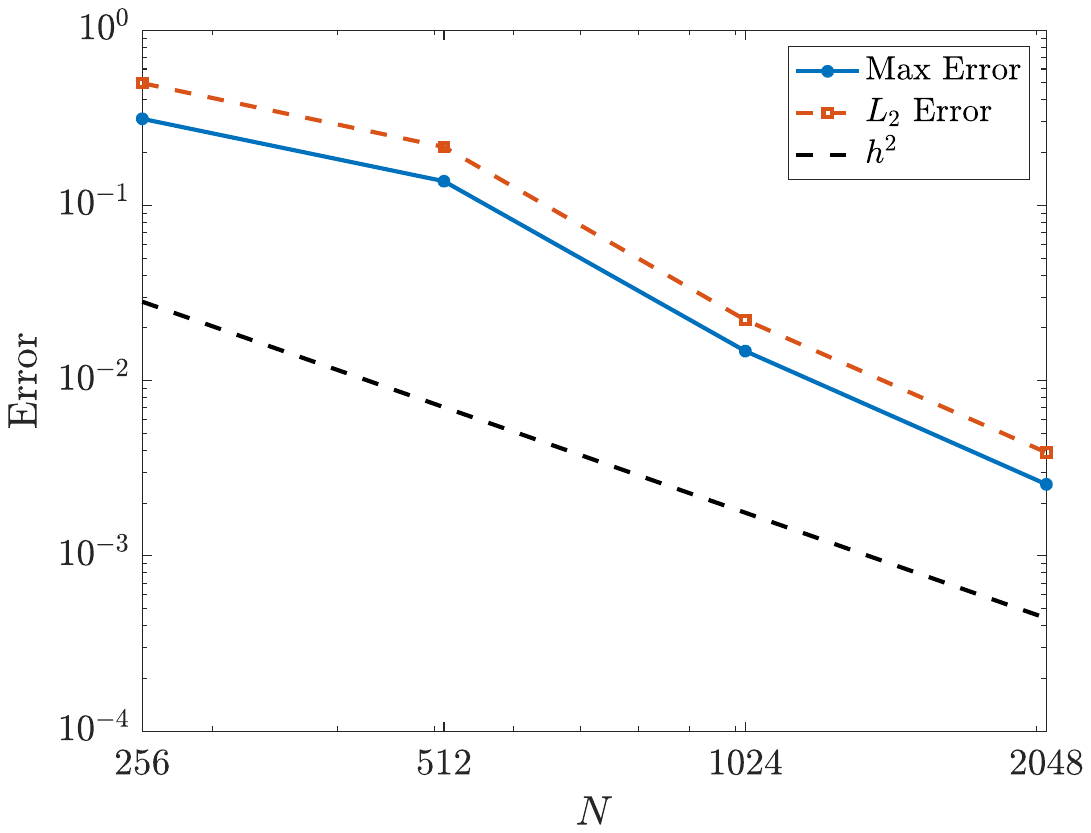}
        \caption{TE}
    \end{subfigure}
    \caption{Self convergence for $u^{\rm sc}$ with two rods}\label{fig:two_rods_conv}
\end{figure}

Fig.~\ref{fig:two_rods_conv} demonstrates the self convergence. The max and $L_2$ errors are computed against the reference solution obtained on mesh with $N=4096$.

\section{Conclusion}
We developed a novel and geometrically robust finite-difference based unfitted lattice Green's function method for exterior scattering in 2D. With the help of local basis functions and lattice Green's functions, the simple finite difference method is enabled to handle irregular geometries elegantly. Our method is analogous to the boundary integral/element method, but does not suffer from singularity in singular integrals. Different boundary conditions are handled in the same numerical framework, and the accuracy and robustness of the method are demonstrated  via a few test examples. 

The developed method can be extended to 3D in a straightforward way -- once we get the lattice Green's function. We did not touch many other aspects of the numerical solutions of Helmholtz equation in this work. Of particular interest is the large wavenumber or frequency, where high order versions will be more relevant. A viable approach is to find the high order lattice Green's function and employ higher degree Lagrange basis functions. Spurious oscillations, corner singularity, etc. fall into other future work. We hope the algorithm in this work serves as a starting point for its development and applications in engineering.



\section*{Acknowledgment}

Q. Xia thanks J. Lai,  S. Tsynkov, and T. Yin for kind and helpful discussions.



\end{document}